\documentclass{amsart}
\usepackage{amssymb, amsthm, txfonts}

\usepackage{graphicx}%
\usepackage{multirow}%
\usepackage{amsmath,amssymb,amsfonts}%
\usepackage{mathrsfs}%
\usepackage[title]{appendix}%
\usepackage{xcolor}%
\usepackage{textcomp}%
\usepackage{manyfoot}%
\usepackage{booktabs}%
\usepackage{algorithm}%
\usepackage{algorithmicx}%
\usepackage{algpseudocode}%
\usepackage{listings}%

\usepackage{bm}
\usepackage{geometry}
\usepackage{tikz}
\usepackage{float}
\usetikzlibrary{intersections, calc, arrows.meta,through,backgrounds}

\theoremstyle{thmstyleone}%
\newtheorem{theorem}{Theorem}
\newtheorem{proposition}[theorem]{Proposition}%
\newtheorem{lemma}[theorem]{Lemma}%

\theoremstyle{thmstyletwo}%
\newtheorem{example}{Example}%
\newtheorem{remark}{Remark}%

\theoremstyle{thmstylethree}%
\newtheorem{definition}{Definition}%

\DeclareMathOperator{\rank}{rank}
\DeclareMathOperator{\sgn}{sgn}

\author{Shuzo Izumi}
\email{sizmsizm@gmail.com}

\newcommand{\BB}[1]{\mathbb{#1}}

\newcommand{\FRAK}[1]{\mathfrak{#1}}
\newcommand{\VEC}{\overrightarrow}

\title[Rotation of a Simplex  
within Another One]{Rotation of a Simplex  
within Another One}

\begin{document}
\maketitle
\begin{abstract}

We are motivated by the most naive problem 
whether a triangle can be rotated within another one. 
We can prove that $\sigma$ can be always moved 
by congruence in an elementary way. 
Generalising this we treat movement of a 
polytope within another one of the same dimension 
in $\BB{R}^n$. By \emph{movement}, we mean 
a non-trivial continuous operation of 
one parameter family of the congruences. 
If we restrict movement to translation, 
the conclusion is very simple. General movements 
are very complicated. The following case is most 
interesting and fairly general. 

Let $\sigma\subset\tau\subset\BB{R}^n$ be 
simplices in Euclidean space such that the 
following holds:\\
$(*)$ All the vertices of $\sigma$ are 
contained in the interior of the different sides of 
$\tau$.

Assume further that $n$ is even. 
Then for a general infinitesimal rotation 
$\boldsymbol\omega$, there exist 
a non-empty open polyhedron $D\subset\BB{R}^n$ 
such that either $\boldsymbol\omega$ or 
$-\boldsymbol\omega$ 
induces a rotation of $\sigma$ within $\tau$ 
with centre in $D$. However, in dimension 3, 
such an infinitesimal rotation is 
not general. Even in these elementary phenomena, 
the parity of the dimension yields difference. 
This is a consequence of the fact that a general 
element of the Lie algebra of the rotation group of 
$\BB{R}^n$ has the full rank if and only if  
$n$ is even. 
\end{abstract} 
\setcounter{section}{-1}
\section{\bf Introduction.}
Take a polytope contained in another one 
of the same dimension. 
We consider whether the former can be moved 
continuously within the latter a little, 
of course not deforming both. 
In this paper, a \emph{polytope} 
means the convex hull of a finite point set in an 
Euclidean space $\BB{R}^n$. 
A \emph{polyhedron} is a subset of 
$\BB{R}^n$ expressed as the intersection of 
a finite number of half spaces bounded by 
affine hyperplanes. 
Then a polytope is just a bounded polyhedron 
(see \cite{bib2}, p.4). We assume as usual 
that $n$-dimensional simplices are always 
non-degenerate, 
that is, their $n+1$ vertices are affinely independent, 
or they are in general position and 
that a polytope and its faces are always closed, 
that is, they contain all of their boundary points 
besides the interior points unless the contrary 
is stated. 

Let $E^+(n)$ denote the group of orientation 
preserving congruences (isometries) of $\BB{R}^n$. 
Let $\sigma\subset\tau\subset\BB{R}^n$ be 
polytopes. By a \emph{movement} of $\sigma$ 
within $\tau$ by $E^+(n)$, 
we mean a short smooth non-trivial continuous curve 
$C:[0,\epsilon)\longrightarrow E^+(n)$ such that 
$C(0)$ is the identity of $E^+(n)$ and  
$C(t)(\sigma)\subset\tau\ (t\in(0,\epsilon))$. 
Here, $C$ need not be a part of a one-parameter 
subgroup of $E^+(n)$.  
This $C$ gives an isotopy of $\sigma$ within $\tau$. 

If $C$ is a subset of a one-parameter subgroup 
of translations (parallel displacement), we call it 
\emph{movement by translation}. In this case 
all points of $\sigma$ move on parallel lines. 
Similarly 
we define \emph{movement by rotation} 
when $C$ is a subset of a 
one parameter subgroup of the special orthogonal 
group $SO(n)$. In this case, the centre is fixed. 
We treat easy facts about movement of polytopes 
in \S \ref{translation}. 

We are mainly concerned with the most interesting 
case: 

\noindent$(*)$ 
Both $\sigma\subset\tau\subset\BB{R}^n$ are 
$n$-dimensional simplices such that 
each of the vertices of $\sigma$ are contained 
in different facets of $\tau$. 

The special orthogonal group SO$(n)$ of dimension $n$ 
is generated by elements of the Lie algebra 
$\FRAK{so}(n)$. We call an element 
$\boldsymbol\omega\in\FRAK{so}(n)$  
an \emph{infinitesimal rotation}. 
We call $\boldsymbol\omega$ as \emph{direction} 
identifying it with its positive constant multiples 
$c\cdot\boldsymbol\omega$ $(c>0)$. 
Note that the axis of an infinitesimal rotation may 
not be unique in dimension greater than 3. 
Our main result is the following: \\\vskip.2ex
{\bf Theorem \ref{main}.} 
\emph{Suppose that $n$ is even and 
$\sigma\subset\tau\subset\BB{R}^n$ 
are simplices with $\dim\sigma=\dim\tau=n$ 
which satisfies the condition $(*)$ above. Then, for a 
general direction $\boldsymbol\omega$, there 
exists an $n$-dimensional open polyhedron 
$D\subset\BB{R}^n$ such that, for any choice of the 
centre in $D$, a movement by rotation of 
$\sigma$ within $\tau$ either in the direction 
$\boldsymbol\omega$ or 
$-\boldsymbol\omega$ is possible.}\vskip2ex

An exceptional impossible case, 
Example \ref{nonadm}, 
will be helpful to understand the meaning of our 
conditions. In this example, rotation is impossible 
but even a small perturbation of the inner simplex 
under the condition $(*)$ makes it possible. 
We cannot extend our arguments to dimension 3 
by Example \ref{3dim}. This difference is brought 
by the fact that the matrices expressing infinitesimal 
rotations have generally the full rank for even $n$ 
in contrast to odd dimensional case. 
\begin{remark}
It is obvious that two inclusions of congruent 
simplices $\sigma$ and $\sigma'$ within another 
$\tau$ are not always isotopic 
by congruence even if the embeddings are quite 
loose and even if they have the same orientation. 
That is, the embedding cannot be always 
continuously deformed to each other. 
Thus the adjective ``small'' is necessary 
for the movement. 
\end{remark}
Supposing the condition $(*)$, our results 
on small movements of a simplex within 
another one are summarised as follows.
$$\begin{array}{|c||c|c|c|}\hline
\dim\tau&\text{translation}&\text{rotation}
&\text{congruence}\\\hline\hline
2                     &\times
&\scalebox{1.5}{\raisebox{-.3mm}{$\circ$}}
&\circledcirc\\\hline
3&\times&\triangle
& ?
\\\hline
\text{even}\ge 4&\times
&\scalebox{1.5}{\raisebox{-.3mm}{$\circ$}}
& ?
\\\hline
\end{array}$$
Here the meanings of the symbols are as follows:

$\circledcirc$: always possible, \quad
\scalebox{1.5}{$\circ$}: generally possible, 
$\times$: impossible,

$\triangle$: neither the possible case nor the
impossible case is rare.\\
\section{\bf Translation and Congruence.}\label{translation}
Let us confirm the intuitively obvious facts 
about movement. 
\begin{proposition}\label{tight}
Let $\sigma$ and $\tau$ be polytopes with 
$\sigma\subset\tau\subset\BB{R}^n$. 
\begin{enumerate}
\item
If every facet $\tau_i$ of $\tau$ contains at least 
one vertex {\rm P}$_i$ of $\sigma$, then $\sigma$ 
cannot be moved by translation. 
\item
Let $\tau:=[\rm{Q}_1\cdots\rm{Q}_n]$ be a 
$n$-dimensional simplex. Here, 
$\rm{Q}_1,\dots,\rm{Q}_n$ are the vertices of 
$\tau$. If the facet $\tau_1$ of $\tau$ 
opposite to $Q_1$ contains no 
point of $\sigma$, $\sigma$ can be moved by 
translation a little within $\tau$ and by rotation 
also. 
\end{enumerate}
\end{proposition}
\begin{proof}
\noindent\begin{enumerate}
\item
Suppose the converse: 
there is a translation $\varphi$ which moves 
$\sigma$ within $\tau$. 
Let $\bf v$ denote a vector of the direction 
of $\varphi$. 
Take an interior point {\rm Q} of $\sigma$. 
Then there exists 
$i\in[0,n]$ and a point {\rm R}$\in\tau_i$ such that 
$\overset{\longrightarrow}{{\rm QR}}$ 
has the same direction as $\bf v$. 
Let us put {\rm Q}$'=\varphi^{-1}({\rm R})$. 
We can take $\varphi$ very small. Then {\rm Q}$'$ is 
an interior point of $\tau$. Putting 
{\rm P}$:=\varphi({\rm P}_i)$, we have a
parallelogram ${\rm Q}'{\rm P}_i{\rm PR}$. 
The midpoint {\rm M} of {\rm PQ}$'$ coincides with 
that of {\rm P}$_i${\rm R}. If {\rm P}$\in\tau$,  
the M is an interior point of $\tau$ as the midpoint 
PQ$'$ and M belongs to the facet $\tau_i$, 
a contradiction. This proves that 
{\rm P}$=\varphi({\rm P}_i)\notin\tau$, impossibility 
of translation of $\sigma$ within $\tau$. 
\item
Let G be the barycentre of the facet $\tau_1$ 
and put ${\bf v_i}:=\VEC{\rm{Q}_1\rm{Q}_i}$. 
Then ${\bf v}_i$ $(i\neq 1)$ form an affine basis of 
$\BB{R}^n$ and $\bf v:=\VEC{\rm{Q}_1G}$ 
is expressed as 
$$
{\bf v}=c_2{\bf v}_2+\cdots+c_n{\bf v}_n\qquad 
(c_i>0,\ c_2+\cdots+c_n)
.$$
Hence $\bf v$ is an inward vector for every facet 
of $\tau$ excepting $\tau_1$. 
Then all vertices of $\sigma$ on the facets 
of $\tau$ can be translated a little in the direction 
$\bf v$ within $\tau$. Since all other vertices of 
$\sigma$ are interior points of $\tau$, 
they can be translated a little in the direction 
$\bf v$ within $\tau$ also. Hence $\sigma$ 
can be translated a little within $\tau$. 
Since a translation is approximated by a rotation, 
the latter is the desirable rotation. (A more 
precise description can be given using the 
exponential map given below. See also 
Theorem \ref{rotation}.) 
\end{enumerate}
\end{proof}
If assume $(*)$, the movement by congruence 
is always possible in dimension 2  as follows. 
The author does not know whether the higher  
dimensional analogue holds or not. 
\begin{proposition}
Let $\sigma:=[\rm{P}_0\rm{P}_1\rm{P}_2]$ and
$\tau:=[\rm{A}_0\rm{A}_1\rm{A}_2]$ be 
2-dimensional simplices (triangles) in $\BB{R}^2$ 
satisfying the condition $(*)$ in Introduction. 
Then $\sigma$ can 
be moved by congruence within $\tau$ a little. 
\end{proposition}
\begin{proof}
We may assume that P$_0$, P$_1$, P$_2$ are 
contained in A$_1$A$_2$, A$_2$A$_0$, A$_0$A$_1$ 
respectively. 
Take $t\in\BB{R}$ with the small absolute value 
and $\rm{P}_0'$ on $\rm{A}_1{A}_2$ such that 
$t=\rm{P}_0\rm{P}_0'$, where $t>0$ implies 
that $\rm{P}_0'$ is contained in the segment 
$\rm{P}_0A_2$ and $t<0$ implies that $\rm{P}_0'$ 
is contained in the segment $\rm{P}_0A_1$. 
Then take $\rm{P}_1'$ on $\rm{A}_2{A}_0$ such that 
$|\rm{P}_0\rm{P}_1|=|\rm{P}_0'\rm{P}_1'|$ as the 
following figure. 
\begin{figure}[H]
\centering
\begin{tikzpicture}
\draw [line width=.6pt]
(0,0)node[left]{A$_1$}--(4,2)node[above]{A$_2$}
--(9,-.5)node[right]{A$_0$}--cycle;
\draw [line width=.6pt]
(2,1)node[above left]{P${}'_0$}
--(8,0)node[above right]{P$'_1$};
\draw [line width=.6pt]
(1,.5)node[above left]{P$_0$}
--(7,.5)node[above right]{P$_1$};
\draw [dotted, line width=.6pt](1,.5)
--(3,-.18)node[below]{P$_2$}--(7,.5);
\draw [dotted, line width=.6pt](2,1)
--(3.8,-.1)node[below]{P$'_2$}--(8,0);
\fill [white] (3,.5) circle [radius=.2];
\draw (3,.5) node{\small$x$};
\fill [white] (6,.5) circle [radius=.2];
\draw (6,.5) node{\small$y$};
\fill [white] (3.37,.8) circle [radius=.2];
\draw (3.37,.8) node{\small$x'$};
\fill [white] (6.35,.2) circle [radius=.2];
\draw (6.37,.23) node{\small$y'$};
\fill [white] (1.45,.6) circle [radius=.2];
\draw(1.5,.6) node{\small$\alpha$};
\fill [white] (6.5,.6) circle [radius=.18];
\draw(6.5,.6) node{\small$\beta$};
\fill [white] (1.7,.84) circle [radius=.2];
\draw(1.7,.84) node{\small$t$};
\fill [white] (4,.62) circle [radius=.2];
\draw(4,.62) node{\small$\theta$};
\draw (4.68,.486) node[below left]{Q}
--(4.63,-.28) node[below]{S};
\draw(4.63,-.04)--(4.86,-.06)--(4.85,-.28);
\end{tikzpicture}
\caption{movement by congruence}
\end{figure}
\noindent This is possible since $|t|$ is taken 
small. Here, $x,y,x',y'$ are the side lengths and 
$\alpha,\beta,\theta$ are the vertex angles of 
the thin triangles with $0<\alpha,\ \beta<\pi/2$. 
The sign of $\theta$ is determined in accordance 
with that of $t$. 
By the theorem of sine we have 
\begin{gather*}
\frac{x'}{\sin\alpha}
=\frac{x}{\sin (\pi-\theta-\alpha)},
\qquad
\frac{y'}{\sin(\pi-\beta)}
=\frac{y}{\sin (\beta-\theta)}.
\end{gather*}
By our assumption, $x'+y'=x+y$. Then we have 
$$
\frac{\sin (\theta+\alpha)}{\sin\alpha}x'
+\frac{\sin (\beta-\theta)}{\sin\beta}y'=x'+y'
,$$
$$
[\sin\beta\sin(\theta+\alpha)
-\sin\alpha\sin\beta]x'
+[\sin\alpha\sin(\beta-\theta)
-\sin\alpha\sin\beta]y'
=0.
$$
Applying  l'Hôpital's rule, we have
\begin{gather*}
\lim_{\theta\to 0}
\frac{y'}{x'}
=-\lim_{\theta\to 0}
\frac{\sin\beta\sin(\theta+\alpha)
-\sin\alpha\sin\beta}{\sin\alpha\sin(\beta-\theta)
-\sin\alpha\sin\beta}
=
\frac{\tan\beta}{\tan\alpha}.
\end{gather*}
Let Q be the point on $\rm{P}_0\rm{P}_1$ 
such that 
$|\rm{P}_0Q|:|Q\rm{P}_1|=\tan\alpha:\tan\beta$ 
and let S be the foot of the perpendicular from 
Q to $\rm{A}_0\rm{A}_1$. 

Suppose that $\rm{P}_2$ is included in the segment 
$\rm{A}_1\rm{S}$. Take P$_2'$ such that 
$\triangle\rm{P_0}\rm{P_1}\rm{P_2}$ is congruent to 
$\triangle\rm{P_0}'\rm{P_1}'\rm{P_2}'$ in this 
order of vertices. Then $\rm{P}_2'$ stay inside 
of $\tau$ for a small positive value of $\theta$. 
Then $\sigma$ can be moved to 
$\triangle\rm{P}'_0\rm{P}'_1\rm{P}'_2$ within $\tau$ 
by the congruent transformation group. 

If $\rm{P}_2$ is included in the segment 
$\rm{A}_2\rm{S}$, $\sigma$ can be moved  
to $\triangle\rm{P}'_0\rm{P}'_1\rm{P}'_2$ 
for a negative value of $\theta$ similarly. 
In particular, 
if $\rm{P}_2=\rm{S}$, $\sigma$ can be moved 
in both directions by congruent transformations 
respectively. 
\end{proof}
\section{\bf Infinitesimal Rotation.}
If an infinitesimal rotation 
$\boldsymbol\omega\in \FRAK{so}(n)$ is given, 
it is known to be an $n$-dimensional linear 
transformation expressed by a real skew 
symmetric matrix $S=S_{\boldsymbol\omega}$ 
of size $n\times n$. If $n$ is even, it is 
well-known that there exists an orthogonal matrix $T$ 
such that  
\[
T^{-1}ST=
\begin{pmatrix}
A_1 & 0 & \dots & 0\\
0 &\ddots &0 &\vdots\\
\vdots&0 &\ddots&0\\
0 &\dots &0&A_{n/2}
\end{pmatrix},
\text{ where } A_i=\begin{pmatrix} 0&-\lambda_i\\
\lambda_i &0\end{pmatrix},\ 
\lambda_i\in\BB{R}.
\]
If $n$ is odd, we have only to add final row and 
final column which totally consist of entity 0. 
Thus $\rank S$ takes the maximum generally and 
the value is $n$ if $n$ is even, and $n-1$ if $n$ 
is odd. 

Let us consider an $n\times n$ matrix $S$ to 
operate on a vertical $n$-vector by a matrix product 
from the left. Recall that 
\[
\exp S
:=E+\frac{S}{1!}
+\frac{S^2}{2!}
+\frac{S^3}{3!}+\cdots
\quad
(E:\text{ the identity matrix})
\]
converges in the normed vector space 
$SO(n)\cong\BB{R}^{n^2}$ 
and defines the exponential map 
\[
\exp: \FRAK{so}(n)\longrightarrow SO(n),
\qquad
(\boldsymbol\omega\text{ represented by }S)
\longmapsto \exp S
.\]
In order to give a rotation on $\BB{R}^n$, 
we need to  
specify the centre Q besides $\exp S$. 
Let $\bf x$ and ${\bf q}$ be the coordinates of 
point X and Q respectively. We treat these 
as vertical vectors in calculations below. 
Then the map of the \emph{rotation} 
$\exp_{\rm Q} tS
=\exp_{\rm Q}S_{t\boldsymbol\omega}$ 
with centre Q generated by 
$t\boldsymbol{\omega}$ $(t\in\BB{R}\setminus\{0\})$ 
is defined as 
\[
\begin{array}{cccl}
\hspace*{.15cm}\BB{R}^n & 
\xrightarrow{\exp_{\rm Q}S_{t\boldsymbol\omega}}
&\hspace{-3cm}\BB{R}^n 
\\
\rotatebox{90}{$\in$}&&\hspace{-3.2cm}\rotatebox{90}{$\in$}
\\
{\bf x}=({\bf x}-{\bf q})+{\bf q} &
\longmapsto & (\exp_{\rm Q} tS)({\bf x}):=
(\exp tS)({\bf x}-{\bf q})+{\bf q}
.\end{array}
\]
This image is approximated by ${\bf x}+tS({\bf x}-{\bf q})$ 
neglecting $t^k$ $(k\ge 2)$ terms. 
Thus the signed moving distance of P in the direction of 
$\bf n$ by $\exp_{\rm Q} tS$ is approximately
\[
{\bf n}\cdot \bigl(tS({\bf x}-{\bf q})\bigr)
=t(S{\bf n})\cdot ({\bf q}-{\bf x})
\]
for $t$ with $|t|<<1$. Here the centre dot 
``$\cdot$'' indicates the inner product of 
vectors. Then we have the following. 
\begin{lemma}\label{rotatable}
Let $\alpha$ be an affine hyperplane of $\BB{R}^n$ 
and take a point {\rm X}$=({\bf x})\in\alpha$ and 
{\rm Q}$=(\bf q)\in\BB{R}^n$. 
Let ${\bf n}$ be a unit normal vector of $\alpha$. 
Let $S$ be a real skew symmetric 
matrix of size $n\times n$. 
Then we have the following. 
\begin{figure}[ht]
\centering
\begin{tikzpicture}
 \coordinate (A) at (0,1.8);
 \coordinate (B) at (4,1.8);%
 \coordinate (C) at (5,.2);%
 \coordinate (D) at (1,.2);%
 \draw [thick] (C)--(D)--(A)--(B)--cycle;
 \coordinate[label=above:${\bf n}$] (n) at (2.5,2.8);
 \filldraw [white] (2.45,1.8) circle [radius=.16cm];     
 \coordinate[label=below right:$\rm X$] (X) at (2.5,1);
 \filldraw [white] (1.25,1.8) circle [radius=.15cm];     
 \filldraw [white] (2.85,1.8) circle [radius=.12cm];     
 \coordinate[label=above:{\rm Q}] (Q) at (.5,2.3);%
 \draw [->, line width=1pt, dashed] (X)--(n);
 \draw [line width=1pt] (Q)--(X);
 \coordinate [label=:$\alpha$] (al) at (4.3,.2);
 \coordinate[label=above right:$S({\bf x}-{\bf q})$] (s) 
 at (3.2,2.5); 
 \draw [->, line width=1pt] (X)--(s);
\end{tikzpicture}
\caption{infinitesimal rotation and a plain}
\end{figure}
\begin{enumerate}
\item
If $S{\bf n}\cdot ({\bf q}-{\bf x})>0$, namely 
if Q is in $\bf n$ side with respect to $\alpha$, 
then 
$(\exp_{\rm Q} tS)({\bf x})$ belongs to $\bf n$ side 
with respect to $\alpha$ for $0<t<<1$. 
\item
If $S{\bf n}\cdot ({\bf q}-{\bf x})<0$, namely 
if Q is in $-\bf n$ side with respect to $\alpha$, 
then 
$(\exp_{\rm Q} tS)({\bf x})$ belongs to $-\bf n$ 
side with respect to $\alpha$ for $0<t<<1$.  
\item
Suppose that $S{\bf n}\cdot ({\bf q}-{\bf x})=0$, 
$S\neq 0$. 
\begin{enumerate}
\item
If ${\bf n}\cdot ({\bf q}-{\bf x})>0$ $($resp. $<0)$, 
then $(\exp_{\rm Q} tS)({\bf x})$ belongs to 
$\bf n$ side $($resp. $-\bf n$ side) for $t$ 
with $0<|t|<<1$. 
\item
If ${\bf n}\cdot ({\bf q}-{\bf x})=0$, then 
${\bf x}\in\alpha$ and $(\exp_{\rm Q} tS)({\bf x})$ 
stays in $\alpha$. 
\end{enumerate}
\end{enumerate}
\end{lemma}
For the proof of the statement (3), (a) above, 
we use the fact that the rotation 
$\exp_{\rm Q} tS\in SO(n)$ preserves 
the distance from Q. 
\section{\bf Rotation of a polytope within 
another one.}%
\label{main section}
\begin{definition}\label{admissible}
Let $\sigma$ and $\tau$ be polytopes such that 
$\sigma\subset\tau\subset\BB{R}^n$, 
$\boldsymbol{\omega}\in\FRAK{so}(n)$ 
an infinitesimal rotation expressed by $S$ and 
a point Q$\in\BB{R}^n$. 
We say that the quadruple 
$(\sigma,\tau,\boldsymbol{\omega},{\rm Q})$ 
\emph{admits rotation}, or it is \emph{admissible} 
if, for any $\bf{x}\in\sigma$, 
$(\exp_Q tS)(\bf{x})\in\tau$ 
for sufficiently small $t>0$.   
\end{definition}
Considering convexity, it is essential for admissibility 
of a polytope to consider the behaviour of vertices of 
$\sigma$ on the boundary $\partial\tau$ of $\tau$.  
Lemma \ref{rotatable} implies the following. 
\begin{theorem}\label{polytope}
Let $\sigma$ be a polytope and 
$\tau$ another polytope with facets 
$\tau_0,\tau_1,\dots,\tau_k$ such that 
$\sigma\subset\tau\subset\BB{R}^n$. 
Let ${\rm P}_0,{\rm P}_1,\dots,{\rm P}_m$ 
be the vertices of $\sigma$ lying on the 
boundary $\partial\tau
:=\tau_0\cup\tau_1\cup\dots\cup\tau_n$ of 
$\tau$. 
Let ${\bf n}^j$ denote the unit inward 
normal vector of the facet $\tau_j$. 
$($Note that one P$_i$ may belongs to many different 
$\tau_j$ and one $\tau_j$ may contains many 
different P$_i$.$)$ 
Given an infinitesimal rotation 
$\boldsymbol{\omega}\neq 0$ expressed by a 
skew symmetric matrix $S$, 
we consider the affine function 
\[
f_{ij}({\bf x})
:=(S{\bf n}^j)\cdot({\bf x}-{\bf p}_i)\qquad
\]
for each $(i,j)$ with ${\rm P}_i\in\tau_j$, 
where ${\bf p}_i$ denotes the coordinate expression 
of {\rm P}$_i$. Then 
$(\sigma,\tau,\boldsymbol{\omega},{\rm Q})$ 
is admissible if and only if either of the following two 
conditions are satisfied for each $(i,j)$ with 
${\rm P}_i\in\tau_j$: 
\begin{enumerate}
\item
$f_{ij}({\bf q})> 0$.
\item
$f_{ij}({\bf q})=0$,\quad
${\bf n}^j\cdot({\bf q}-{\bf p}_i)\ge 0$.
\end{enumerate}
\end{theorem}
\begin{theorem}\label{rotation}
Let $\sigma$ and 
$\tau$ be polytopes such that 
$\sigma\subset\tau\subset\BB{R}^n$. 
Suppose that 
the facets $\tau_h$, $\tau_{h+1}$, \dots, $\tau_m$ 
$(1\le h\le n)$ contain no vertices of $\sigma$. 
Let ${\bf n}^i$ denote the unit internal 
normal of $\tau_i$. Take an infinitesimal rotation 
$\boldsymbol{\omega}$ which is expressed by 
a skew symmetric matrix $S$. Suppose that 
$S{\bf n}^0,S{\bf n}^1,\dots,S{\bf n}^{h-1}$ are 
linearly independent. Then there exist 
two non-bounded open polyhedrons 
$\Omega,\, \Omega'\subset\BB{R}^n$ such 
that the quadruple 
$(\sigma,\tau,\boldsymbol{\omega},{\rm Q})$ 
(resp. $(\sigma,\tau,-\boldsymbol{\omega},{\rm Q}')$) 
is admissible for any Q$\in\Omega$ 
(resp. Q$'\in\Omega'$). 
\end{theorem}
\begin{proof}
Let 
P$_{i1}$, P$_{i2}$,\dots, P$_{ik_i}$ be all the vertices 
of $\sigma$ on $\tau_i$ $(0\le i\le h-1)$. 
Of course, 
P$_{ij}=$P$_{i'j'}$ may happen for some $i\neq i'$. 
Only these P$_{ij}$ concern small rotation of $\sigma$. 
Take the affine functions 
\[
f_{ij}({\bf x}):=(S{\bf n}^i)\cdot({\bf x}-{\bf p}_{ij})
\qquad (0\le i\le h-1,\ 1\le j\le k_i).
\]
Here, ${\bf p}_{ij}$ is the coordinate expression 
of P$_{ij}$. Then, for every fixed $i\in[0, h-1]$, 
the differences of $f_{ij}({\bf x})=0$ $(1\le j\le k_i)$ 
are constants and define mutually parallel hyperplanes. 
The functions 
\[
(y_1,y_2,\dots,y_h):=(f_{01},f_{11},\dots,f_{h-1\ 1})
\]
form \emph{a part of} affine coordinates system 
of $\BB{R}^n$ by the linear independence 
of $S{\bf n}^0,S{\bf n}^1,\dots,S{\bf n}^{h-1}$. 
The coordinates $y_1,y_2,\dots,y_h$ can take 
arbitrary values independently. Hence, by giving 
$y_1,y_2,\dots,y_h$ sufficiently large values, 
all $f_{ij}$ take arbitrarily large positive values 
simultaneously. This means that 
\[
\bigcap_{i,j}D_{ij}^+\neq \emptyset\quad 
(D_{ij}^+:=\{{\bf x}:\sgn f_{ij}({\bf x})=+\})
\]
is a non-bounded open polyhedron i.e. an intersection 
of half spaces. 
For any Q$\in \bigcap_{i,j}D_{ij}^+$, the quadruple 
$(\sigma,\tau,\boldsymbol{\omega},{\rm Q})$ is 
admissible by Lemma \ref{rotatable}. 
Quite similarly, we have 
$\bigcap_{ij}D_{ij}^-\neq\emptyset$ and 
$(\sigma,\tau,-\boldsymbol{\omega},{\rm Q'})$ 
is admissible for any Q$'\in \bigcap_{i,j}D_{ij}^-$. 
\end{proof}
\section{\bf An Example.}
We show a critical example of inclusion of 
simplices $\sigma\subset\tau\subset\BB{R}^n$ 
which admits no movement but a slight displacement 
of $\sigma$ makes it admissible for a suitable 
choice of the centre. This will explain 
the conditions of Theorem \ref{main}.
\begin{example}\label{nonadm}
Let $\sigma:=[\rm{P}_0\rm{P}_1\rm{P}_2]$ and
$\tau:=[\rm{A}_0\rm{A}_1\rm{A}_2]$ be 
2-dimensional simpleces such that vertices of 
$\sigma$ are interior points of different 
facets of $\tau$ as in Fig.\,2. 
Let ${\bf n}_i^0$ denote the inward normals 
of the facets of $\tau$ at $P_i$ $(i=0,1,2)$. 
If all ${\bf n}_0^0,{\bf n}_1^0,{\bf n}_2^0$ meet 
at one point R of $\tau$, we can rotate 
$\sigma$ within $\tau$ in both 
direction, clockwise and counter clockwise. 

Next, we assume that all 
${\bf n}_0^0,{\bf n}_1^0,{\bf n}_2^0$ meet 
at one point R outside $\tau$ as Fig.\,3, 
where $+$ sign put near the normal line P$_0$R 
implies that, if the centre Q is taken in this side, 
${\rm P}_0$ admits clock-wise small rotation 
$\boldsymbol{\omega}$ inside with 
respect to the side ${\rm A}_1{\rm A}_2$.  
The $-$ sign implies that ${\rm P}_0$ goes out 
of the side ${\rm A}_1{\rm A}_2$. 
Similar for other normal lines. 
We associate to every point of $\BB{R}^2$ 
the sequence of its three signs determined by its 
position relative to these normals. 
\begin{figure}[h]
\centering
\begin{tikzpicture}
\draw [thick] (-.35,-.5)  arc (30:150:1.4);
\draw [thick] (.4,0)  arc (30:150:2.4);
\draw [thick] (.77,0)  arc (30:150:2.8 8);
\draw [thick]
 (-4.3,.33) node[left] {A$_1$}-- 
(-2.2,1.7) node[above] {A$_2$}-- 
(1.4,.05) node[right]{A$_0$}--cycle;

\filldraw [black] (-3.15,1.07) circle [radius=.07cm];  
\filldraw [black] (-.6,.95) circle [radius=.07cm];  
\filldraw [black] (-1.5,.2) circle [radius=.07cm];   
\filldraw [black] (-1.6,-1) circle [radius=.07cm];   

\coordinate[label=left:R] (R) at (-1.6,-1);
\coordinate[label=below:\, \, P$_1$] (P1) at (-.6,.95);
\coordinate[label=left:P$_0\,$] (P0) at (-3.15,1.07);%
\coordinate[label=below left:\, P$_2$] (P2) at (-1.5,.2);%
\draw [line width=1pt] (P0)--(P1)--(P2)--cycle; 

\coordinate (d) at ($(R)!3.5cm!(P1)$);
\coordinate[label=left:{${\bf n}_1^0$}] (e) at ($(P1)!3cm!(R)$);
\draw [thick, dashed] (d)--(e);
\draw [thick, ->] (P1)--(e);

\coordinate (f) at ($(R)!4cm!(P0)$);
\coordinate[label=right:{${\bf n}_0^0$}] (g) 
at ($(P0)!3.5cm!(R)$);
\draw [thick, dashed] (f)--(g);
\draw [thick, ->] (P0)--(g);

\coordinate[label=right:{${\bf n}_2^0$}] (h) at ($(R)!3.5cm!(P2)$);
\coordinate (k) at ($(P2)!2cm!(R)$);
\draw [thick, dashed] (P2)--(k);
\draw [thick, ->] (P2)--(h);

\coordinate[label=left:+] (s0) at (-3.6,1.6);
\coordinate[label=left:$-$] (t0) at (-3.15,1.7);

\coordinate[label=left:+] (s1) at (-.25,1.6);
\coordinate[label=left:$-$] (t1) at (.2,1.6);

\coordinate[label=left:+] (s2) at (-1,2);
\coordinate[label=left:$-$] (t2) at (-1.3,2);
\end{tikzpicture}
\caption{seldom non-rotatable case in dimension 2}
\end{figure}
Then there is no point with simple sign 
sequence $(+++)$ nor $(---)$ with respect to 
$f_0,f_1,f_2$. As a sole point with simply 
non-negative or simply non-positive sign sequence, 
R has the sign sequence $(0\ 0\ 0)$. 
But R is not contained in $\tau$. 
This implies that 
$(\sigma,\tau,\boldsymbol{\omega},{\rm R})$ 
and $(\sigma,\tau,-\boldsymbol{\omega},{\rm R})$ 
are not admissible. Indeed, P$_2$ goes out of 
$\tau$ even by a small rotation with centre R. 

If P$_2$ is taken a little lefter on A$_1$A$_2$, 
there appears a little triangle bounded by 
the normal lines of sides of $\tau$, which 
consists of points with sign sequence $(+++)$.
If we take Q in that triangle, 
$(\tau,\sigma,\boldsymbol{\omega},\rm{Q})$ 
is admissible for clockwise $\boldsymbol{\omega}$. 
If P$_2$ is taken a little righter on A$_1$A$_2$, 
there appears a little triangle of points with sign 
sequence $(---)$. If we take Q in that triangle, 
$(\tau,\sigma,-\boldsymbol{\omega},\rm{Q})$ 
is admissible. 
\end{example}
\section{\bf Rotation of an even dimensional simplex 
in another one.}\label{mainsection}
Let $\sigma:=[{\rm P}_0{\rm P}_1\cdots {\rm P}_n]
\subset\BB{R}^n$ 
be an $n$-dimensional simplex $(n\ge 2)$, where 
${\rm P}_0,{\rm P}_1,\dots, {\rm P}_n$ denotes 
the vertices. 
The facet opposite to P$_i$ is expressed as 
\[
\sigma_i:=[{\rm P}_0{\rm P}_1\cdots 
\check{{\rm P}}_i\cdots {\rm P}_n]\qquad
(i=0,1,\dots,n).
,\]
where the check mark ``\ $\check{ }$\ '' implies 
omission. Thus $\sigma_i$ is an $(n-1)$-dimensional 
simplex. We need the following well-known fact 
that the unit inward unit normals 
of its facets ${\bf n}^0,{\bf n}^1,\dots,{\bf n}^n$
are in general position. 
Indeed, this follows from the fact that the 
intersection of the subspace parallel to 
$n$ facets 
$\sigma_0,\sigma_1,\dots,\check{\sigma_i},
\dots,\sigma_n$ is $0$-dimensional. 

Let $\sigma\subset\tau\subset\BB{R}^n$ be 
$n$-dimensional simplices. 
Suppose that $\boldsymbol{\omega}$ is an 
infinitesimal rotation expressed by a full rank skew 
symmetric matrix $S$. Then 
$S{\bf n}^0,S{\bf n}^1,\dots,S{\bf n}^n$ 
are in general position as well as the inward normals 
${\bf n}^0,{\bf n}^1,\dots,{\bf n}^n$ of facets 
of $\sigma$. 
Excluding the non-general case when some 
vertex of $\sigma$ is included in a face of $\tau$ of 
codimension greater than 1, we may assume the 
condition $(*)$ in view of Theorem \ref{rotation}. 
Then, we have the following, our main 
result. 
\begin{theorem}\label{main}
Let $n$ be an even natural number. 
Take $n$-dimensional simplices 
$\tau:=[{\rm A}_0{\rm A}_1\cdots {\rm A}_n]$ and 
$\sigma:=[{\rm P}_0{\rm P}_1\cdots {\rm P}_n]$ 
with $\sigma\subset\tau\subset\BB{R}^n$. 
We suppose the condition $(*)$ in Introduction: 
all vertices ${\rm P}_i$ 
of $\sigma$ are contained in the interior 
$\overset{\circ}{\tau_i}$ of the different 
facets $\tau_i$ of $\tau$. Let 
$\boldsymbol{\omega}$ be a direction with  
the full rank matrix expression $S$. We put 
\[
f_i({\bf x}):=(S{\bf n}^i)\cdot({\bf x}-{\bf p}_i),\quad
D_i^0:=f_i^{-1}(0)\subset\BB{R}^n
\qquad(0\le i\le n).
\]
where ${\bf p}_j$ denotes the coordinate expression 
of {\rm P}$_j$. Then, for a general choice 
of ${\bf p}_i$, $D_i^0$ do not meet at a point. 
In that case at least one of the open polyhedrons 
\[
D_{[0,n]}^+:=
D_0^+\cap D_1^+\cap\dots\cap D_n^+,\quad
D_{[0,n]}^-:=
D_0^-\cap D_1^-\cap\dots\cap D_n^-
\]
is non-empty and we have the following.
\begin{enumerate}
\item
If $D_{[0,n]}^+\neq\emptyset$, then for any 
$\mathrm{Q}\in D_{[0,n]}^+\cup\big(\overline{D_{[0,n]}^+}
\cap\tau\big)$, the quadruple 
$(\tau,\sigma,\boldsymbol{\omega},\rm{Q})$
is admissible. 
\item
If $D_{[0,n]}^-\neq\emptyset$, then for any 
$\mathrm{Q}\in D_{[0,n]}^-\cup\big(\overline{D_{[0,n]}^-}
\cap\tau\big)$, the quadruple 
$(\tau,\sigma,-\boldsymbol{\omega},\rm{Q})$
is admissible. 
\end{enumerate}
\end{theorem}
\begin{proof}
Since $S{\bf n}^1,S{\bf n}^2,\dots,S{\bf n}^n$ 
are linearly independent as stated above, 
the functions 
\[
(y_1,y_2,\dots,y_n):=(f_1,f_2,\dots,f_n)
\]
form an affine coordinate system of $\BB{R}^n$. 
The coordinate hyperplanes $D_1^0,D_2^0,\dots,D_n^0$ 
for $(y_1,y_2,\dots,y_n)$ are transversal and meet at 
one point, say R. 
Since the coordinates $y_1,y_2,\dots,y_n$ 
are independent, there are two non-empty 
polyhedrons 
\[
D_{[1,n]}^+:=
D_1^+\cap D_2^+\cap\dots\cap D_n^+,
\quad
D_{[1,n]}^-:=
D_1^-\cap D_2^-\cap\dots\cap D_n^-
\]
with the simple sign sequence 
$(++\dots+)$ and $(--\dots-)$ respectively. 

The hyperplane $D_0^0$ can be expressed as 
\[
f_0({\bf y})
:=a_1y_1+a_2y_2+\dots+a_ny_n+a=0\quad
((a_1,a_2,\dots,a_n)\neq(0,0,\dots,0))
\]
with respect to the affine coordinates 
$(y_1,y_2,\dots,y_n)$. 
Since $D_i^0$ $(1\le i\le n)$ do not meet 
at one point, we see that $a\neq 0$. 
\begin{enumerate}
\item
Suppose that $a_i>0$ $(i\neq 0)$ for certain $i$. 
Fixing the values of $y_j$ 
$(j\neq i)$ positive (resp. negative) and making 
$y_i$ very large (resp. small), 
we have a point of $D_{[0,n]}^+$ (resp. $D_{[0,n]}^-$). 
Thus $D_{[0,n]}^+\neq\emptyset$ 
and $D_{[0,n]}^-\neq\emptyset$. 
\item
Suppose that all $a_i\le 0$ $(i\neq 0)$. 
Since $D_0^0\neq\emptyset$, 
there exists $a_j<0$ $(j\neq 0)$. 
\begin{enumerate}
\item
If $a>0$ further, 
taking $y_i$ $(i\neq 0)$ with sufficiently small 
positive values, we have $f_0({\bf y})>0$. Hence 
$D_{[0,n]}^+\neq\emptyset$. 
If $y_i<0$ for all $i\neq 0$ we have $f_0({\bf y})>a$. 
Hence $D_{[0,n]}^-=\emptyset$. 
\item
If $a<0$, taking negative $y_i$ $(i\neq 0)$ with 
sufficiently 
small absolute values, we have $f_0({\bf y})<0$. Hence 
$D_{[0,n]}^-\neq\emptyset$. If $y_i>0$ $(i\neq 0)$, 
we have $f_0({\bf y})\le a<0$ and 
$D_{[0,n]}^+=\emptyset$. 
\end{enumerate}
\end{enumerate}
Hence, at least 
one of $D_{[0,n]}^+$ and $D_{[0,n]}^-$ is non-empty. 
If $Q$ belongs to $D_{[0,n]}^+$ or to 
$\partial D_{[0,n]}^+\cap\tau$, then 
$(\sigma,\tau,\boldsymbol{\omega},{\rm Q})$ 
admits rotation by Lemma \ref{rotatable}. 
If $Q$ belongs to $D_{[0,n]}^-$ or to 
$\partial D_{[0,n]}^-\cap\tau$, then 
$(\sigma,\tau,-\boldsymbol{\omega},{\rm Q})$ 
admits rotation. 

Finally, we confirm that generally all $D_i^0$ 
$(0\le i\le n)$ do not meet at one point. 
To see this exactly, first 
determine the intersection point of $D_i^0$ 
$(1\le i\le n)$ by solving the equation 
\[
(S{\bf n}^i)\cdot({\bf x}-{\bf p}_i)=0\qquad
(1\le i\le n)
\]
with respect to the original euclidean coordinates 
$\bf x$. By independence of 
$S{\bf n}^1,S{\bf n}^2,$ 
$\dots,S{\bf n}^n$, this has a 
unique solution. The components of ${\bf x}$ 
are expressed as linear combinations of the 
components of  ${\bf p}_1,{\bf p}_2,\dots,{\bf p}_n$. 
Substituting $\bf x$ in the remaining equation 
$(S{\bf n}^0)\cdot{\bf x}=(S{\bf n}^0)\cdot{\bf p}_0$ 
by these solutions, we see that all the components of 
${\bf p}_0$ are expressed by components  
of ${\bf p}_1,{\bf p}_2,\dots,{\bf p}_n$. 
This is a non-trivial necessary condition for all 
$D_i^0$ to meet at a single point, because the 
coordinates of point ${\bf p}_0$ is independent of 
those of ${\bf p}_1,{\bf p}_2,\dots,{\bf p}_n$ 
originally. Thus this condition defines a thin 
closed subset of the product space 
$(\partial\tau)^{n+1}$. 
Hence, all $D_i^0$ $(0\le i\le n)$ 
do not meet at a point generally. 
\end{proof}
\section{\bf Rotation of a 3-dimensional simplex 
in another one.}
In the 3-dimensional case, we cannot say that, 
for general given 
$(\tau,\sigma,\boldsymbol{\omega})$ 
satisfying $(*)$ in Introduction, 
there exists $Q$ such that either 
$(\tau,\sigma,\boldsymbol{\omega},\rm{Q})$ or 
$(\tau,\sigma,-\boldsymbol{\omega},\rm{Q})$ 
is admissible as Theorem \ref{main}. 
This is seen in the following example. 
The author guesses that the similar thing happens 
for all odd $n\ge 5$. 
\begin{example}\label{3dim}
There exist 3-dimensional simplices 
$\sigma\subset\tau\subset\BB{R}^3$ such that 
the condition $(*)$ in Introduction is satisfied 
and there is an open set 
$\Omega\subset\FRAK{so}(3)$ of directions 
and for any $\boldsymbol{\omega}\in\Omega$ 
both the quadruples 
$(\tau,\sigma,\boldsymbol{\omega},{\rm Q})$ and 
$(\tau,\sigma,-\boldsymbol{\omega},{\rm Q})$ 
are not admissible for any Q$\in\Omega$.

To see this, take a regular tetrahedron 
[A$_0$A$_1$A$_2$B$_3$] and take a point A$_3$ 
on the extension of its edge 
A$_0$B$_3$ very near to B$_3$. We obtain a 
new simplex $\tau:=[$A$_0$A$_1$A$_2$A$_3]$. 
Take points P$_i$ on the facets 
$\tau_i$ opposite to A$_i$ $(i=0,1,2,3)$ as illustrated 
in Fig.\,4, left. Let $\boldsymbol{\omega}$ be the 
perpendicular vector from A$_0$ to the plane 
$\alpha$ determined by A$_1$,A$_2$,B$_3$. 
It represents the infinitesimal rotation clockwise 
seen from the above. 
Projecting P$_0$, P$_1$, P$_2$, P$_3$, 
A$_0$, A$_3$ to the plane $\alpha$, we obtain 
P$_0'$, P$_1'$, P$_2'$, P$_3'$, A$_0'$, A$_3'$ 
respectively in Fig.\,4, right. Let $S$ be the matrix 
expressing $\boldsymbol{\omega}$ and 
${\bf n}^i$ the unit inward normal of $\tau_i$
and put
\[
f_i({\bf x}):=(S{\bf n}^i)\cdot({\bf x}-{\bf p}_i),\quad
D_i^0:=f_i^{-1}(0)\subset\BB{R}^n
\qquad(0\le i\le 4).
\]
Since the hyperplane $D_i^0$ includes ${\bf n}^i$ 
and $\boldsymbol{\omega}$, 
they are projected as dot-dash lines 
on $\alpha$. The meanings of the signs 
close to the dot-dash lines are as follows: 

\begin{figure}[H]

\begin{tikzpicture}
\draw [line width=1pt]
(2,3.5) node[right] {A$_0$}-- 
(4,1) node[right] {A$_2$}-- 
(2,0) node[below]{A$_1$}--
(-.36,.5) node[left] {A$_3$}--cycle;
\draw [line width=1pt] (4,1)--(-.36,.5);
\coordinate [label=left:B$_3$] (B_3) at (0,1);
\filldraw [white] (.75,.6) circle [radius=.1cm]; 
\draw [line width=1pt,dashed] (2,0)--(B_3)--(4,1);
\filldraw [white] (2,1) circle [radius=.1cm]; 
\filldraw [white] (2,.8) circle [radius=.1cm]; 
\draw [line width=1pt] (2,0)--(2,3.5);

\draw [line width=1pt, black] (2.2,.6) circle [radius=.07cm]; 
\coordinate [label=right:P$_0$] (P0) at (2.2,.6); 
\draw [line width=1pt, black] (1.8,2.4) circle [radius=.07cm]; 
\coordinate [label=left:P$_1$] (P1) at (1.8,2.4); 
\filldraw [black] (1.2,1.3) circle [radius=.07cm]; 
\coordinate [label=above:P$_2$] (P2) at (1.2,1.3); 
\filldraw [black] (2.3,1.85) circle [radius=.07cm]; 
\coordinate [label=right:P$_3$] (P3) at (2.3,1.85); 

\filldraw [black] (.3,-1) circle [radius=.07cm]; 
\coordinate [label=right:this side] (vis) at (.3,-1); 
\draw [line width=1pt, black] (2.5,-1) circle [radius=.07cm]; 
\coordinate [label=right:far side] (invis) at (2.5,-1); 

\draw [arrows=->,line width=1pt] (.1,3.4)--(.1,1.5);

\draw (-.18,2.7) node {$\omega$};

\end{tikzpicture}
\begin{tikzpicture}
\draw [dashed,line width=1pt]
(3.8,3) node[right] {A$_2$}-- 
(-.3,3.3) node[left] {A$_3'$}--
(2,0) node[below]{A$_1$};

\draw [line width=1pt]
(3.8,3)-- (2,0);

\draw [line width=1pt] (2,0)--(2,2)--(3.8,3);
\draw [line width=1pt] (2,2)--(-.3,3.3);
\draw (2,2.25) node [below right] {A$_0'$};

\draw [line width=1pt] (.27,.03) arc (340:20:.5);
\draw [line width=1pt] (.09,.39)--(.29,.31)--(.34,.46);
\draw (-.2,.2) node {$\omega$};

\coordinate [label=above:P$_0'$] (P0) at (2.2,2.5); 
\draw [line width=1pt,dash dot] (.3,3.5)--(3.9,1.6);
\draw (3.9,1.85) node {$+$};
\draw (3.7,1.55) node {$-$};
\coordinate [label=:$D_0^0$] (D0) at (4.2,1.1);
\filldraw [white] (2.2,2.5) circle [radius=.07cm]; 
\draw [line width=1pt, black] (2.2,2.5) circle [radius=.07cm]; 

\filldraw [black] (1.7,2.5) circle [radius=.07cm]; 
\coordinate [label=left:P$_1'$] (P1) at (1.7,2.5); 
\draw [line width=1pt,dash dot] (1.7,0)--(1.7,3.6);
\draw (1.5,3.4) node {$+$};
\draw (1.87,3.4) node {$-$};
\coordinate [label=:$D_1^0$] (D1) at (1.7,3.5);

\filldraw [black] (1.2,1.86) circle [radius=.07cm]; 
\coordinate [label=below:P$_2'$] (P2) at (1.2,1.86); 
\draw [line width=1pt,dash dot] (3.8,3.2)--(.1,1.3);
\draw (.5,1.27) node {$+$};
\draw (.3,1.6) node {$-$};
\coordinate [label=:$D_2^0$] (D2) at (-.18,.9);

\filldraw [black] (2.3,1.5) circle [radius=.07cm]; 
\coordinate [label=below:P$_3'$] (P3) at (2.3,1.5); 
\draw [line width=1pt,dash dot] (3.4,.9)--(-.1,2.8);
\draw (3.26,1.2) node {$+$};
\draw (3.06,.9) node {$-$};
\coordinate [label=:$D_3^0$] (D3) at (3.7,.45);
\coordinate [label=:\phantom{a}] (b) at (2,-1);
\end{tikzpicture}
\caption{non-rotatable case in dimension 3}
\end{figure}
\vskip-5ex

The $+$ sign means that, if the centre Q of the 
rotation is taken in this $+$ side with respect to 
$D_i^0$, a small rotation 
induced by $\boldsymbol{\omega}$ preserves P$_i$ 
inside $\tau$ with respect to $\tau_i$. 
The $-$ sign means that, if the centre Q of the 
rotation is taken in this $-$ side with respect to 
$D_i^0$, P$_i$ leaves $\tau$ through $\tau_i$ 
by a small rotation 
generated by $\boldsymbol{\omega}$. 

There are no point where all the signs are 
non-negative 
and no point where all signs are non-positive. 
Hence there is no admissible rotation of the 
simplex $[\rm{P}_0\rm{P}_1\rm{P}_2\rm{P}_3]$ for any 
choice of the centre Q.  
This situation cannot be changed by a 
sufficiently small perturbation of 
$(\sigma,\boldsymbol{\omega})$. 
Hence, this pair is an interior point 
(see Remark \ref{semianalytic}) of the set 
of points which permit no admissible quadruples 
by any choice of the centre Q. 
On the other hand, the set of rotatable quadruple 
have also interior points. Indeed, if we move 
the point P$_2$ to a position very near 
to A$_3$ on the facet $\tau_2$, the region of points 
with sign sequence $(++++)$ appears and for any 
Q in that region we have an admissible quadruple 
$(\tau,\sigma,\boldsymbol{\omega},{\rm Q})$. 
\end{example}
\begin{remark}\label{semianalytic}
For fixed $\tau$, we can naturally define 
a semianalytic structure \cite{bib1} 
on the set of the pairs 
$(\sigma,\boldsymbol{\omega})$ for a fixed $\tau$ in 
Example \ref{3dim}, which bears a topology. 
\end{remark}

\end{document}